\pdfoutput=1
\documentclass[pdf, 10pt]{article}
\usepackage[square,authoryear]{natbib}
\usepackage{marsden_article}
\usepackage{amscd}

\input xy
\xyoption{all}

\def\E{\mathbb{E}}

\begin{document}

\title{Stochastic Variational Integrators}

\author{Nawaf Bou-Rabee\thanks{Applied \& Computational Mathematics (ACM), Caltech, Pasadena, CA 91125 ({\tt nawaf@acm.caltech.edu}).}  \and 
Houman Owhadi\thanks{Applied \& Computational Mathematics (ACM) and Control \& Dynamical Systems (CDS), Caltech, Pasadena, CA 91125 (\tt{owhadi@acm.caltech.edu}).} }

\maketitle

\begin{abstract}
This paper presents a continuous and discrete Lagrangian theory for stochastic Hamiltonian systems on manifolds.  The main result is to derive stochastic governing equations for such systems from a critical point  of a stochastic action.   Using this result the paper derives Langevin-type equations for constrained mechanical systems and implements a stochastic analog of Lagrangian reduction.  These are easy consequences of the fact that the stochastic action is intrinsically defined.   Stochastic variational integrators (SVIs) are developed using a discretized stochastic variational principle.  The paper shows that the discrete flow of an SVI is a.s.~symplectic and in the presence of symmetry a.s.~momentum-map preserving. A first-order mean-square convergent SVI for mechanical systems on Lie groups is introduced.  As an application of the theory, SVIs are exhibited for multiple, randomly forced and torqued rigid-bodies interacting via a potential.
\end{abstract}

 \section{Introduction}

Since the foundational work of Bismut [1981], the field of stochastic geometric mechanics is emerging in response to the demand for tools to analyze the structure of continuous and discrete mechanical systems with uncertainty \citep{Bi1981, Li1997, LiWa2005, MiReTr2002, MiReTr2003, LaOr2007a, LaOr2007b,MaWi2007, CiLeVa2008}.   Our specific goal within this context is to develop efficient, structure-preserving integrators for long-time simulations of stochastic Hamiltonian systems on manifolds.   Our strategy is to extend variational integrators to this class of systems.

\paragraph{Variational Integrators}

Variational integration theory derives integrators for mechanical systems from discrete variational principles \citep{Ve1988,Ma1992,WeMa1997,MaWe2001}.   The theory includes discrete analogs of the Lagrangian, Noether's theorem, the Euler-Lagrange equations, and the Legendre transform.  Variational integrators can readily incorporate holonomic constraints (via Lagrange multipliers) and non-conservative effects (via their virtual work)  \citep{WeMa1997, MaWe2001}.    Altogether, this description of mechanics stands as a self-contained theory of mechanics akin to Hamiltonian, Lagrangian or Newtonian mechanics.

\begin{figure}[ht!]
\begin{center}
\includegraphics[scale=0.25,angle=0]{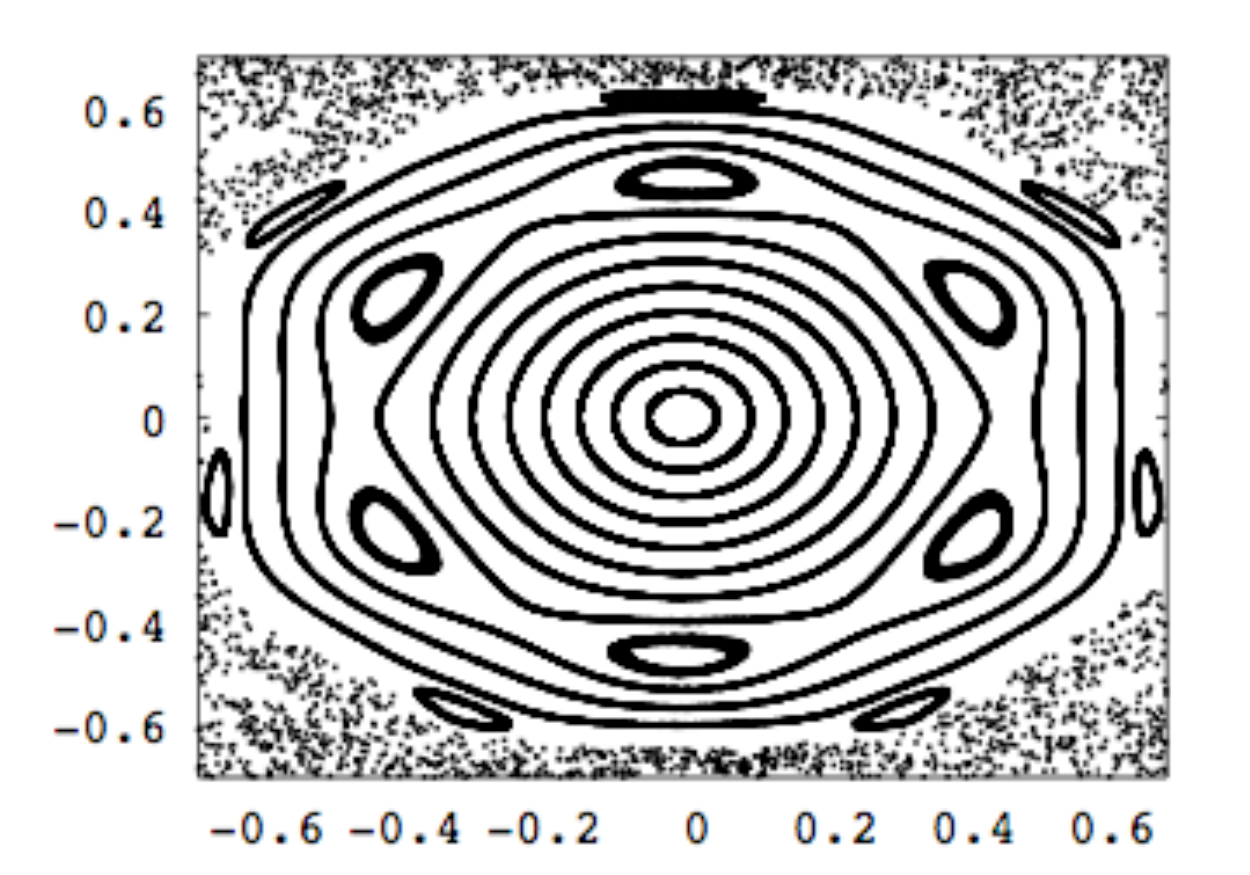}
\includegraphics[scale=0.25,angle=0]{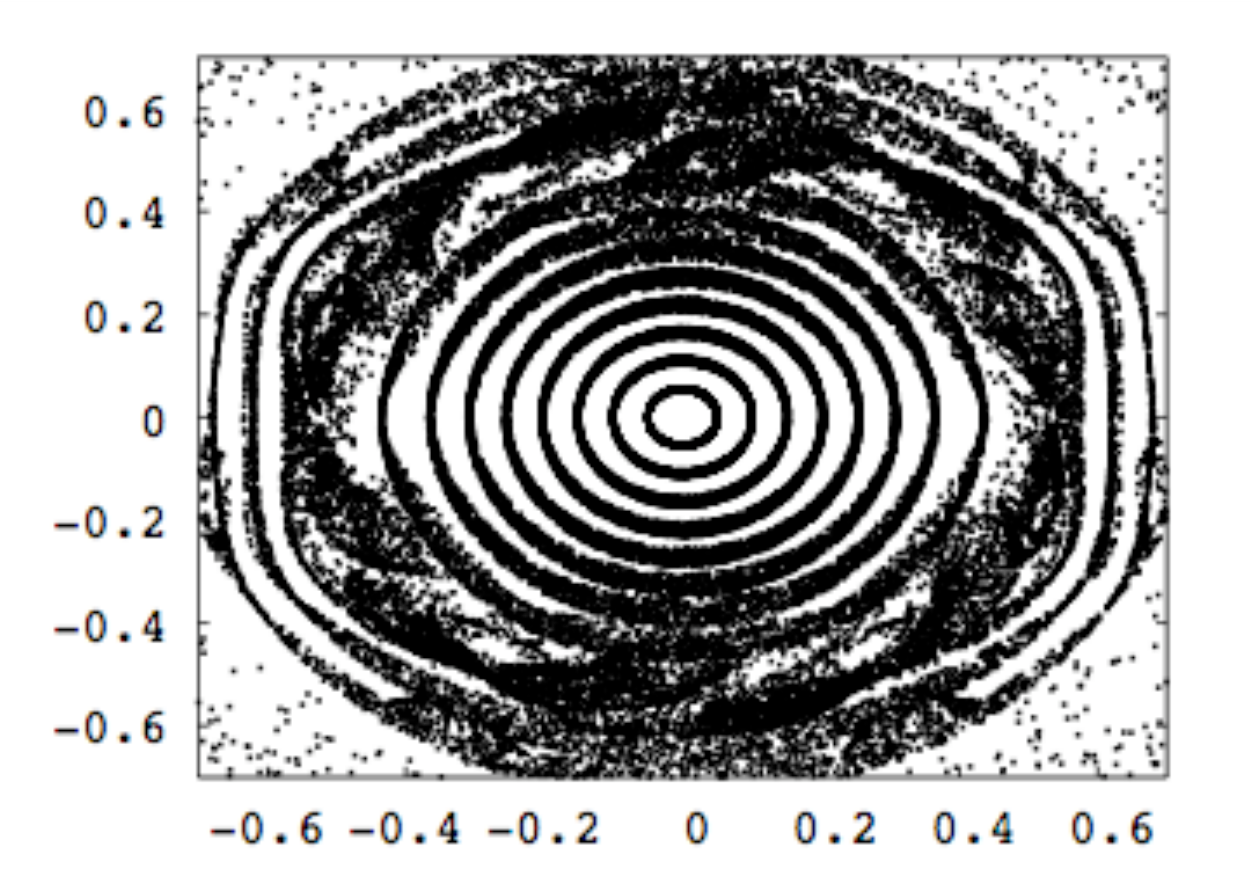} 
\includegraphics[scale=0.25,angle=0]{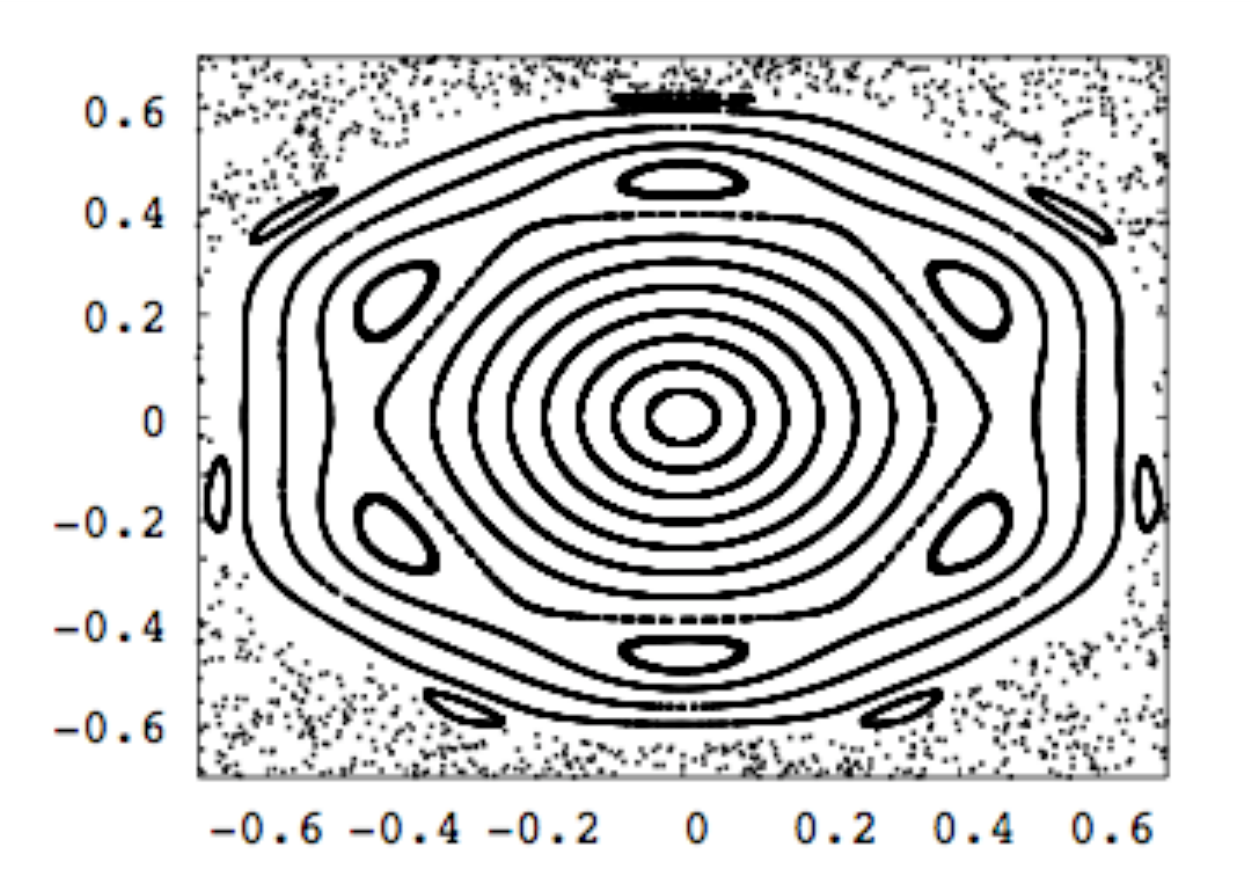}
\includegraphics[scale=0.25,angle=0]{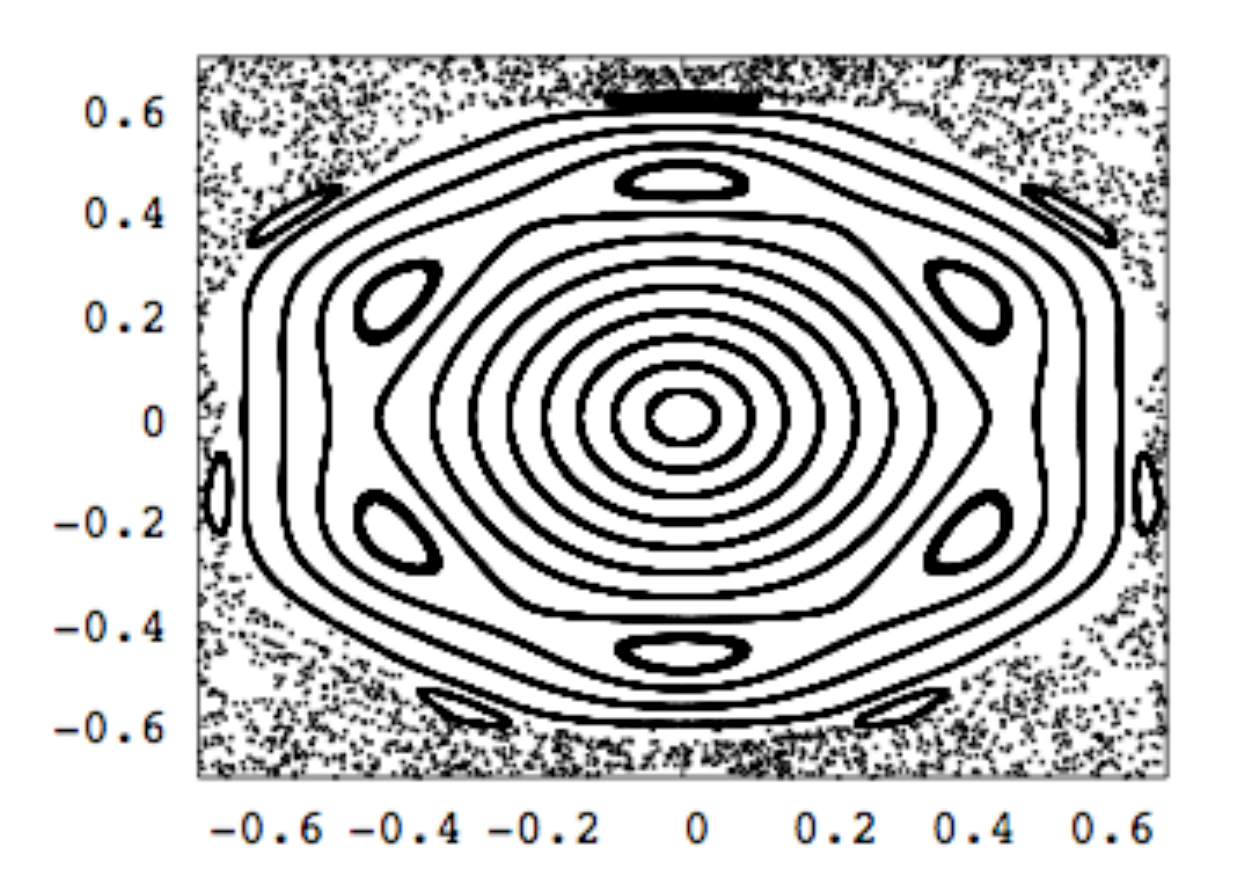} 
\hbox{\hspace{0.0in}   $\begin{array}{c} \text{(a) RK4} \\ $h=0.025$ \end{array}$ \hspace{0.6in}  $\begin{array}{c} \text{(b) RK4} \\ $h=0.05$ \end{array}$ \hspace{0.6in}  $\begin{array}{c} \text{(c) VE} \\ $h=0.025$ \end{array}$ \hspace{0.6in}  $\begin{array}{c} \text{(d) VE} \\ $h=0.05$ \end{array}$ } 
\caption{ \small {\bf Underwater Vehicle Dynamics \citep{BoMa2007}.} 
This figure shows a computation of Poincar\'{e} sections using a second-order accurate variational Euler integrator (VE)  as compared to fourth-order accurate Runge-Kutta (RK4).  Both methods agree with the benchmark at the finer stepsize $h=0.025$.  However, at the coarser stepsize $h=0.05$, RK4 corrupts chaotic invariant sets while the lower-order accurate VE method preserves the structure of the benchmark.    } \label{fig:uv} \end{center}
\end{figure}

One of the distinguishing features of variational integrators is their ability to compute statistical properties of mechanical systems, such as in computing Poincar\'e sections, the instantaneous temperature of a system, etc.  For example, as a consequence of their variational construction, variational integrators are symplectic \citep{Vo1956,Ru1983,Fe1986}.   A single-step integrator applied to a mechanical system is called {\em symplectic} if the discrete flow map it defines exactly preserves the continuous symplectic 2-form and is otherwise called standard.  Using backward error analysis one can show that symplectic integrators applied to Hamiltonian systems nearly preserve the energy of the continuous mechanical system for exponentially long periods of time and that the modified equations are also Hamiltonian (for detailed exposition see \citep{HaLuWa2006}).  Standard integrators often introduce spurious dynamics in long-time simulations, e.g., artificially  corrupt chaotic invariant sets as illustrated in a computation of a Poincar\'{e} section of an underwater vehicle in Fig.~\ref{fig:uv}  using a fourth-order accurate Runge-Kutta (RK4) method and a variational Euler (VE) method \citep{BoMa2007}.


\begin{figure}[ht!]
\begin{center}
\includegraphics[scale=0.75,angle=0]{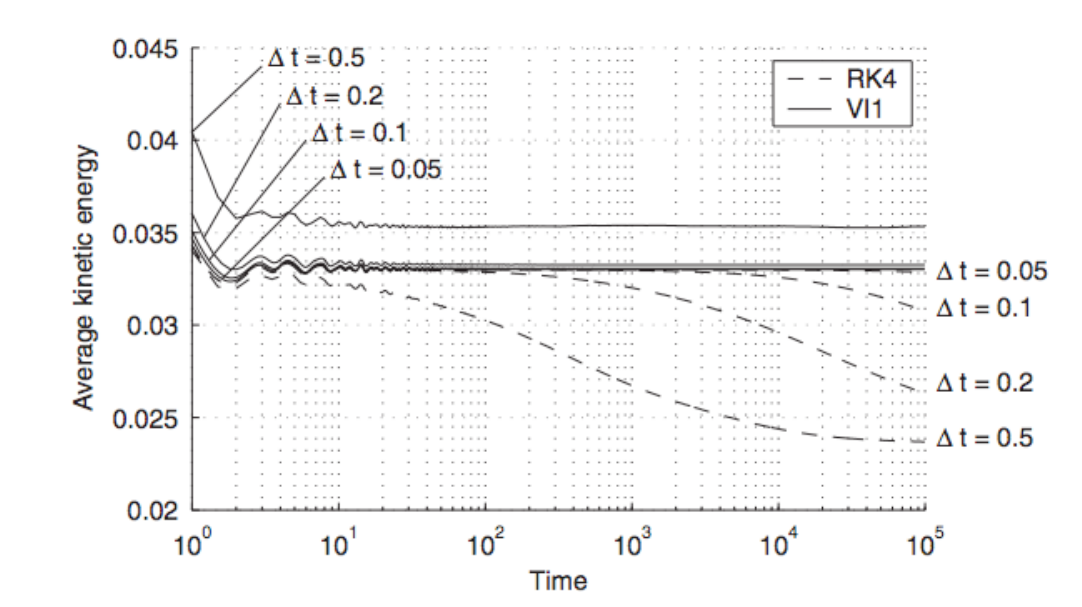}
\caption{ \small {\bf Coupled Spring--Mass Lattice \citep{LeMaOrWe2004b}.} 
This figure shows several plots of the time-averaged instantaneous temperature as a function of time using a first-order accurate variational integrator (VI1) and RK4.    The mechanical system is a coupled spring-mass lattice.  Observe the drift in the computed temperature by RK4 as compared to VI1.
   } \label{fig:springmasslattice} \end{center}
\end{figure}

In addition to correctly computing chaotic invariant sets and long-time excellent energy behavior, evidence is mounting that variational integrators correctly compute other statistical quantities in long-time simulations.  For example, in a simulation of a coupled spring-mass lattice, Lew et al.~found that variational  integrators correctly compute the time-averaged instantaneous temperature (mean kinetic energy over all particles) over long-time intervals, whereas standard methods (even a higher-order accurate one) exhibit a artificial drift in this statistical quantity.    See Fig.~\ref{fig:springmasslattice} and \citep{LeMaOrWe2004b}.   These structure-preserving properties of variational integrators motivate their extension to stochastic Hamiltonian systems.


\paragraph{Main Results}

In his foundational work, Bismut showed that the stochastic flow of certain randomly perturbed Hamiltonian systems with configuration space $\mathbb{R}^n$ extremize a stochastic action.  He called such systems {\bfi stochastic Hamiltonian systems}, and used this property to prove symplecticity and extend Noether's theorem to such systems \citep{Bi1981}.  Mean-square symplectic integrators for stochastic Hamiltonian systems with configuration space $\mathbb{R}^n$ and driven by Wiener processes have been developed \citep{MiReTr2002, MiReTr2003}.

Bismut's work was further enriched and generalized to manifolds by recent work \citep{LaOr2007a, LaOr2007b}.    Lazaro-Cami and Ortega showed that stochastic Hamiltonian systems on manifolds extremize a stochastic action defined on the space of manifold-valued semimartingales \citep{LaOr2007a}.  Moreover, they performed reduction of stochastic Hamiltonian systems on the cotangent bundle of a Lie group to obtain stochastic Lie-Poisson equations \citep{LaOr2007b}.  However, as far as we can tell the converse to Bismut's original theorem, namely a critical point of a stochastic action satisfies stochastic Hamilton's equations, has not been proven.   In fact, as pointed out by Lazaro-Cami and Ortega, a counterexample can be constructed to prove the converse of Bismut's theorem is not true {\em for a certain choice of stochastic action} \citep{LaOr2007a}.

In this paper we restrict our attention to stochastic Hamiltonian systems driven by Wiener processes and assume that the space of admissible curves in configuration space are of class $C^1$.   From the viewpoint of randomly perturbed mechanical systems, this latter restriction is reasonable since random effects often do not appear in the kinematic equation, but rather the balance of momentum equation as white noise forces and torques.  It should be mentioned that the ideas in this paper can be readily extended to stochastic Hamiltonian systems driven by more general semimartingales, but for the sake of clarity we restrict to Wiener processes.  Within this context the results of the paper are as follows:

\begin{itemize}
\item For a class of mechanical systems whose configuration space is a paracompact manifold and which is subjected to multiplicative white noise forces and torques, the paper proves a.s.~a curve satisfies stochastic Hamilton's equations if and only if it extremizes a stochastic action.   This theorem is the main result of the paper.
\item The paper derives governing SDEs for stochastic Hamiltonian systems with holonomic constraints using a constrained stochastic variational principle and for stochastic Hamiltonian systems with nonconservative effects in the drift using a Lagrange-d'Alembert-type principle (for deterministic treatments see \citep{MaRa1999}). The paper performs Lagrangian reduction for stochastic Hamiltonian systems whose configuration space is a Lie group and provides stochastic Euler-Poincar\'e/Lie-Poisson equations for such systems (for deterministic treatment see \citep{MaRa1999}).   These are easy consequences of the fact that the stochastic action is intrinsically defined.  
\item The paper shows how to discretize variational principles to obtain stochastic variational integrators (SVIs), stochastic RATTLE-type integrators for constrained stochastic Hamiltonian systems, and stochastic Euler-Poincar\'e/Lie-Poisson integrators for stochastic Hamiltonian systems on Lie Groups (for deterministic treatments see \citep{MoVe1991,WeMa1997, MaPeSh1998, HaLuWa2006}).  In addition, the paper describes how to derive quasi-symplectic methods for rigid-body-type systems at uniform temperature.
\end{itemize}

\paragraph{Organization of the Paper}

Sufficient conditions for existence, uniqueness, and mean-squared differentiability of stochastic flows on manifolds are recalled in \S \ref{sec:stochasticflows}. In \S \ref{sec:stochastichpmechanics} we extend the Hamilton-Pontryagin (HP) principle to the stochastic  setting to prove a class of mechanical systems with multiplicative noise appearing as forces and torques possess a variational structure.  It should be emphasized (and it is explained in the section) the mechanical system could evolve on a nonlinear configuration space and involve holonomic constraints or nonconservative effects in the drift.  The HP viewpoint is adopted since it unifies the Hamiltonian and Lagrangian descriptions of the system.  By left-trivializing this principle, we also show how to perform Lagrangian reduction in this stochastic setting for stochastic rigid-body-type systems.

In \S \ref{sec:stochasticvi}, SVIs are derived from an abstract discrete Lagrangian and the structure of the resulting discrete flow map is analyzed.  In \S \ref{sec:stochasticveuler}, we concretely show how to design a single-step, stochastic variational Euler integrator for mechanical systems whose configuration space is a Lie group using a simple stochastic discrete HP principle.   If the configuration space is the translation group, the governing SDEs and integrators are in one-to-one correspondence with so-called stochastic Hamiltonian systems and mean-square symplectic integrators  \citep{Bi1981, MiReTr2002, MiReTr2003};  and with the addition of dissipation in the drift term, in one-to-one correspondence with Langevin equations and so-called quasi-symplectic methods  for such systems.   These mean-square symplectic integrators have been numerically tested and shown to possess excellent properties for long-time simulations of mechanical systems governed by Langevin-type equations; see \citep{MiTr2004}.

\begin{figure}[ht!]
\begin{center}
\includegraphics[scale=0.35,angle=0]{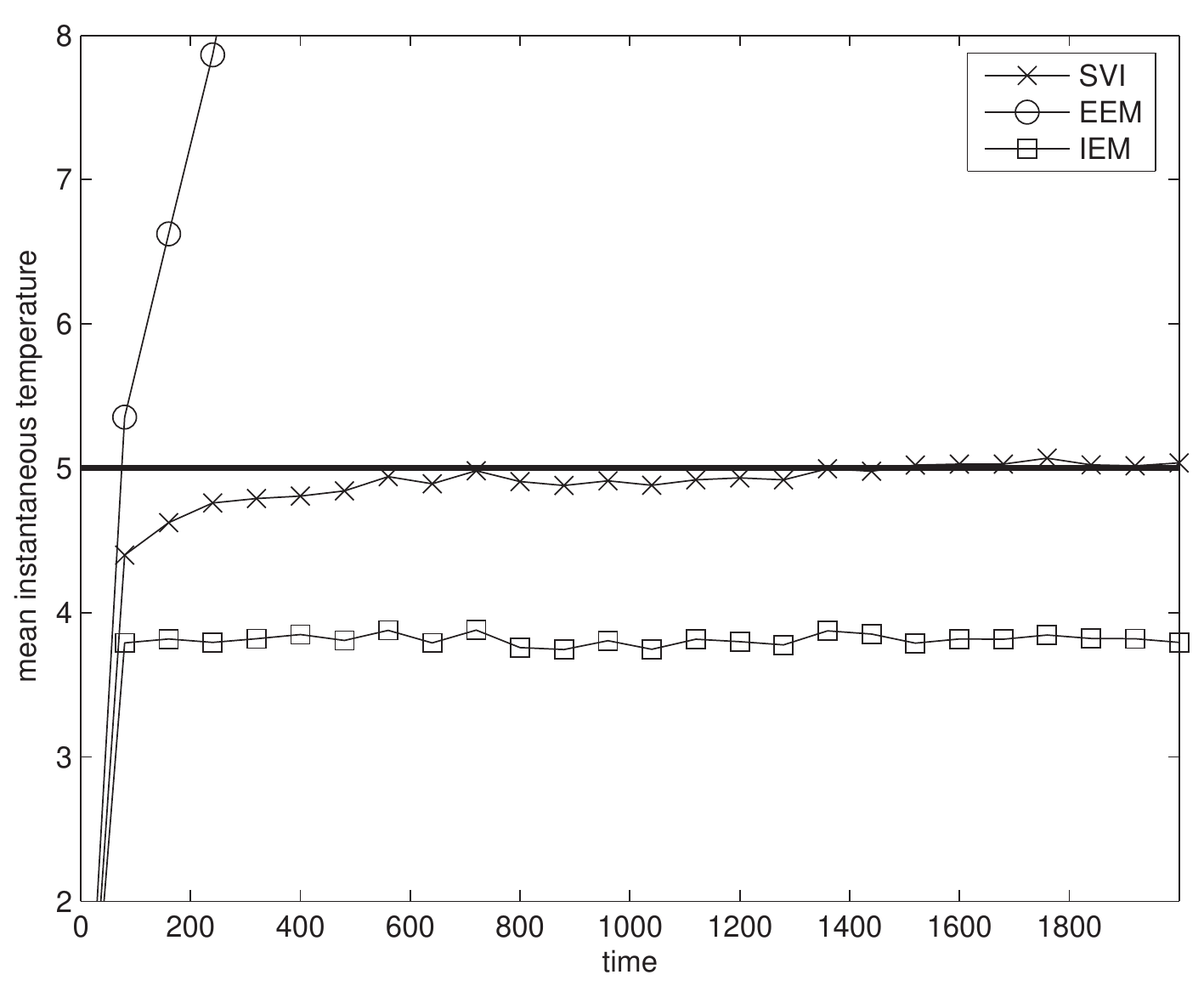}  
\includegraphics[scale=0.35,angle=0]{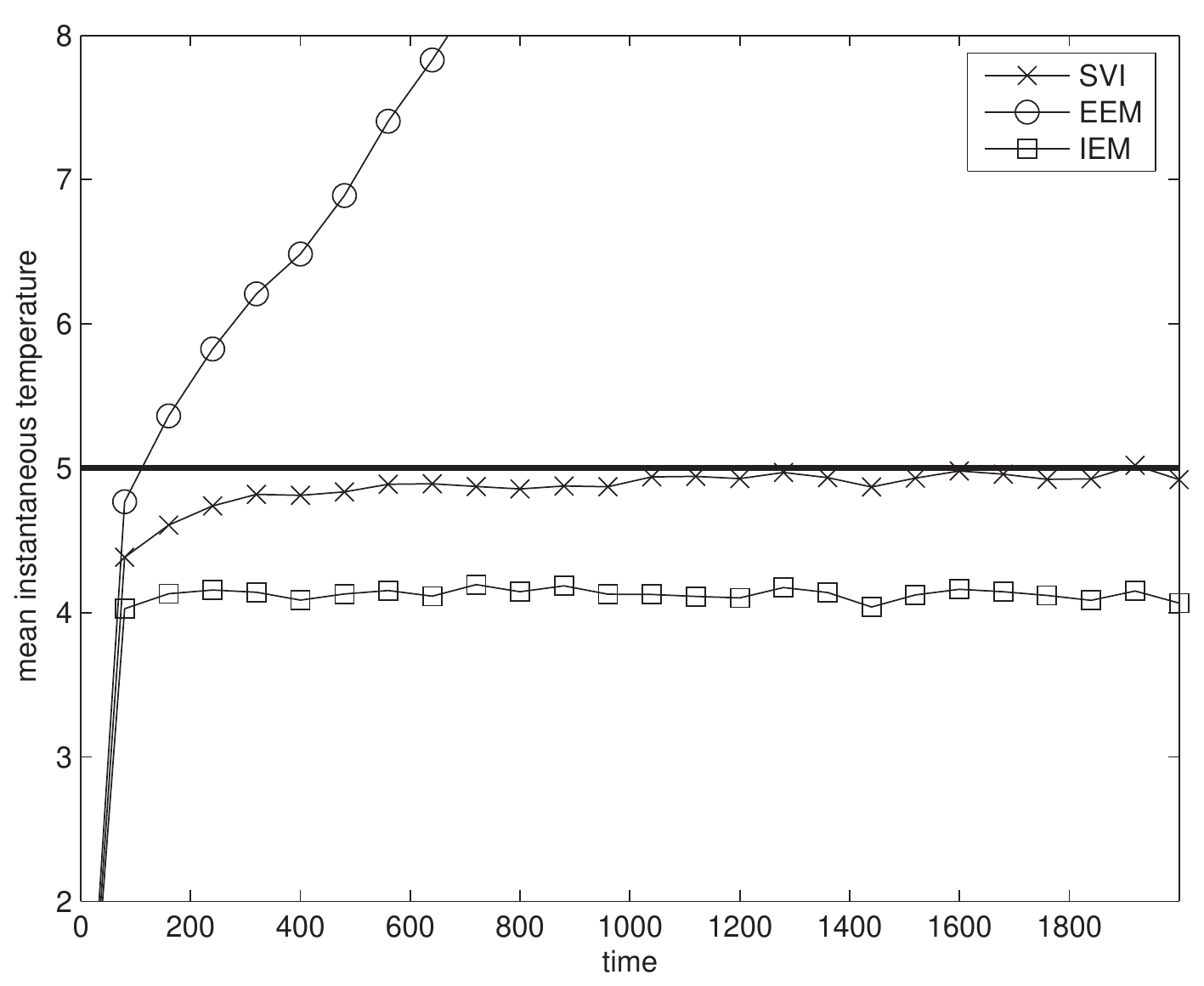}  \\
\hbox{\hspace{1.5in}   (a) $h=0.1$  \hspace{1.2in}  (b) $h=0.05$ } 
\caption{ \small {\bf Ballistic Pendulum at Uniform Temperature  \citep{BoOw2007c}.}  Plots of the mean instantaneous temperature (kinetic energy) of a ballistic pendulum computed using an SVI, explicit Euler-Maruyama (EEM), and implicit Euler Maruyama (IEM) for time-steps $h$ as indicated.   The correct temperature is indicated by the solid line.  Observe that the EEM and IEM schemes artificially heat and cool the system respectively.  A key feature of the ballistic pendulum is that the diffusion and drift matrices associated to the momentums are degenerate, yet the system is still at uniform temperature.   
   } \label{fig:ballistic} \end{center}
\end{figure}

Our own simulations confirm those findings, and will be discussed in \citep{BoOw2007c, BoOw2007b, BoOwWh2007}.   A sample of such results is provided in Fig.~\ref{fig:ballistic}.  It compares an SVI to standard, presumably non-variational methods on a ballistic pendulum at uniform temperature \citep{BoOw2007c}.  The figure shows that an SVI correctly computes the temperature of the system (defined as the mean of the instantaneous temperature with respect to realizations) whereas explicit and implicit Euler-Maruyama schemes (EEM, IEM) do not.   All of these methods are first order strongly convergent and driven by the same realization of Wiener processes.   This computation suggests that an SVI has favorable long-time energy properties whereas EEM and IEM artificially heat and cool the system respectively.    On the other hand, this paper focuses on SVI theory and the structure-preserving properties of SVIs.

In \S \ref{sec:rigidbodieslangevin}, as an application we explain how one can add multiplicative white noise forces and torques to multiple rigid bodies in a fashion that preserves variational structure.  With the addition of dissipation these become Langevin-type equations.  An SVI is provided for such systems.   Using an Ito-Taylor expansion, it is easy to check that the method is first order in the mean-squared sense \citep{KlPl1992}.


 \section{Stochastic Flows on Manifolds}     \label{sec:stochasticflows}

Some standard results on flows of SDEs on manifolds are reviewed here for the reader's convenience.  The reader is referred to the following textbooks on the subject for more detailed exposition \citep{El1982, Em1989, IkWa1989, Ku1990}.  This section parallels the treatment of deterministic flows on manifolds found in Chapter 4 of \citep{AbMaRa2007}.

We start by introducing notation for deterministic vector fields on manifolds which are an important component of SDEs on manifolds.  Let $M$ be an n-manifold.  Recall, a vector field on $M$ is a section of the tangent bundle $TM$ of $M$.  The set of all $C^k$ vector fields on $M$ is denoted by $\mathfrak{X}^k(M)$.

The notion of a probability space is introduced in order to extend the definition of a dynamical system to incorporate noise.   A stochastic dynamical system consists of a base flow on the probability space which propagates the noise and a stochastic flow on $M$ which depends on the noise.

 \begin{definition}[Stochastic Dynamical System]
 A {\bfi stochastic dynamical system} consists of a base flow on a probability space $(\Omega, \mathcal{F}, \mathbb{P})$ and a stochastic flow on a manifold $M$.  The {\bfi base flow} is a $\mathbb{P}$-preserving map, $\theta: \mathbb{R} \times \Omega \to \Omega$, which satisfies: 
\begin{enumerate}
\item
$\theta_0 = \text{id}_{\Omega}: \Omega \to \Omega$, is the identity on $\Omega$
\item 
for all $s, t \in \mathbb{R}$ the group property: $\theta_s \circ \theta_t = \theta_{s+t}$
\end{enumerate}
Given times $0 \le r \le s \le t $,  the {\bfi stochastic flow} on $M$ is a map $\varphi_{t,s}:  \Omega \times M \to M$ such that
\begin{enumerate}
\item
for a.a. $\omega \in \Omega$, the map $(s, t, \omega, x) \mapsto \varphi_{t,s}(\omega) x$ is continuous in $s$, $t$ and $x$
\item
$\varphi_{s,s}(\omega) = \text{id}_M: M \to M$, is the identity map on $M$ for all $s \in \mathbb{R}$
\item
$\varphi$ satisfies the cocycle property:
\[
\varphi_{t,s}(\theta_s(\omega)) \circ \varphi_{s,r} ( \omega) = \varphi_{t,r}( \omega )
\]
\end{enumerate}
\label{def:stochasticflow}
 \end{definition}


This paper is concerned with stochastic dynamical systems that come from {\bfi stochastic laws of motion}; that is, ones whose stochastic flows define solutions of SDEs.   The Stratonovich definition of stochastic integrals is adopted to extend SDEs to a manifold $M$.  The main advantage of the Stratonovich approach being that the chain rule holds for the Stratonovich differential.   Consider a manifold $M$ modelled on a Banach space $E$, and vector fields $X_i \in \mathfrak{X}^k(M)$ for $i=0,...,m$.  Let  $\mathcal{F}_t$ be a nondecreasing family of $\sigma$-subalgebras of $\mathcal{F}$,  and $(W_i(t, \omega), \mathcal{F}_t)$, $i=1,...,m$, be independent Wiener processes for $0 \le t \le T$,.  In terms of these objects, the Stratonovich SDE the paper considers takes the form,
\begin{equation}
d z = X_0(z) dt + \sum_{i=1}^m X_i(z) \circ dW_i,~~~z(0) = z_0 \text{.}
\label{eq:stratonovichsde}
\end{equation}
$X_0$ is referred to as the {\bfi drift vector field} and $X_i$, $i=1,...,m$, are the {\bfi diffusion vector fields}.  A {\bfi Stratonovich integral curve} of (\ref{eq:stratonovichsde}) is a $C^0$-map $c(\cdot, \omega): [0, T]  \to M$ which satisfies:
\[
c(t, \omega) = z_0 + \int_0^t X_0(c(s, \omega)) ds + \sum_{i=1}^m \int_0^t X_i(c(s,\omega)) \circ d W_i (s, \omega) \text{.}
\]

Uniqueness of solutions to (\ref{eq:stratonovichsde}) will be defined in the pathwise sense.

\begin{definition}[Pathwise Uniqueness]
Let $c$ be a Stratonovich integral curve of (\ref{eq:stratonovichsde}).  {\bfi Pathwise uniqueness} of $c$ means that if $\bar{c}: I \to M$ is also a solution to (\ref{eq:stratonovichsde}) on the same filtered probability space with the same Brownian motion and initial random variable, then 
\[
P \left( c(t, \omega) = \bar{c}(t, \omega), ~~\forall~~  0 \le t \le T \right) =1 \text{.}
\]
\end{definition}

Differentiability of the flow map on $M$ will be defined in the mean-squared sense.  In the following we define mean-squared derivatives on the model space $E$ with the understanding that this notion can be extended to $M$ using a local representative of the flow map.

\begin{definition}[Mean-Squared Derivative]
The mean squared norm of $f: E \times \Omega \to E$ is given by:
\[
\| f (x, \omega) \| = \left( \E \left( |f(x, \omega)|^2 \right) \right)^{1/2}  \text{.}
\]
Using this norm one can define the derivative of $f$ in the standard way, i.e., $f$ is mean squared differentiable at $a \in E$ if there is a bounded linear map $Df(a, \omega): E \to E$ that satisfies,
\[
\lim_{\delta \to 0} \frac{ \left\| f(a+\delta, \omega) - f(a, \omega) - D f(a, \omega) \cdot \delta \right\| }{ \| \delta \| } \to 0  \text{.}
\]
\label{def:meansquared}
\end{definition}

As is standard the explicit dependence on the point $\omega \in \Omega$ will usually be suppressed.  With these definitions one can state the following key, but standard theorem \citep{El1982, Em1989, IkWa1989, Ku1990}.

\begin{theorem}[\textbf{Existence, Uniqueness, and Smoothness}]   Let $M$ be a manifold with model space $E$.  Suppose $X_i \in \mathfrak{X}^k(M)$, $i=0,...,m$ and $k \ge 1$, are uniformly Lipschitz, and measureable with respect to $x \in M$.  Let $I=[0, T]$.  Then,
\begin{enumerate}
\item For each $u \in M$, there is a.s.~a $C^0$-curve $c:I \to M$ such that $c(0) = u$ and $c$ satisfies (\ref{eq:stratonovichsde}) for all $t \in I$.  This curve $c : I \to M$ is called a {\bfi maximal solution}.
\item $c$ is pathwise unique.
\item There is a.s.~a mapping $F: I \times M \to M$ such that the curve $c_u: I \to M$ defined by $c_u(t) = F_t(u)$ is a curve satisfying (\ref{eq:stratonovichsde}) for all $t \in I$.  Moreover, $F$ is $C^k$ in the mean-squared sense in $u$ and $C^0$ in $t$.  
\end{enumerate}
\label{thm:maximalsolution}
\end{theorem}

\section{Stochastic HP Mechanics}  \label{sec:stochastichpmechanics}

In this section a variational principle is introduced for a class of stochastic Hamiltonian systems on manifolds.    The action is simply the classical action stochastically perturbed by stochastic integrals.  The key feature of this principle is that one can recover stochastic Hamilton's equations for these systems.  Roughly speaking, this is accomplished by means of taking variations of this action within the space of curves only (not the probability space), and imposing this ``partial differential'' of the action to be zero.

\paragraph{Stochastic HP Principle}

Let $(\Omega, \mathcal{F}, P)$ be a probability space.  Fix an interval $[a,b] \subset \mathbb{R}$ and let $\mathcal{F}_t$ be a nondecreasing family of $\sigma$-subalgebras of $\mathcal{F}$,  and $(W_i(t), \mathcal{F}_t)$, $i=1,...,m$, be independent Wiener processes for $ t \in [a,b]$.  Consider a mechanical system with Lagrangian $\mathcal{L}: TQ \to \mathbb{R}$ and a paracompact, configuration manifold $Q$.

The paper adopts an HP viewpoint to describe this mechanical system with random perturbations.  The HP principle unifies the Hamiltonian and Lagrangian descriptions of a mechanical system  \citep{YoMa2006a, YoMa2006b, Bo2007, BoMa2007}.  The classical HP action integral will be perturbed using {\bfi stochastic potentials}: $\gamma_i: Q \to \mathbb{R}$ for $i=1,...,m$.  Roughly speaking, in the stochastic context the HP principle states the following critical point condition on $TQ \oplus T^*Q$, 
\[
\delta \int _a^b \left[ \mathcal{L}(q, v) dt + \sum_{i=1}^m \gamma_i (q) \circ d W_i + \left\langle p, \frac{dq}{dt} - v \right\rangle dt \right]  = 0,
\]  
where $(q(t),v(t),p(t)) \in TQ \oplus T^*Q$ are varied arbitrarily and independently with endpoint conditions $q (a) $ and $q (b)$ fixed. This principle builds in a Legendre transform, stochastic Hamilton's equations and stochastic Euler--Lagrange equations.    The action integral in the above principle, consists of two Lebesgue integrals with respect to the Lebesgue measure $dt$ and $m$ Stratonovich stochastic integrals with respect to Wiener processes.  This action is random, i.e., for every sample point $\omega \in \Omega$ one will obtain a different, time-dependent Lagrangian system.   However, each system possesses a variational structure which we will make precise in this section.  For a deterministic treatment of time-dependent continuous and discrete Lagrangian systems the reader is referred to \citep{MaWe2001}.

\begin{definition}
The {\bfi Pontryagin bundle} is defined as the Whitney sum $T Q \oplus T^* Q$.     Fixing the interval $[a,b]$, define {\bfi path space}  as
 \begin{align*}
\mathcal{C}([a,b],q_1,q_2) = \{ (q,v,p) \in C^0( [a,b] , TQ \oplus T^*Q) ~|~ q \in C^1([a,b], Q),~q(a)=q_1,~q(b)=q_2 \}  \text{.}
\end{align*}
Let $\mathfrak{G}: \Omega \times \mathcal{C}([a,b],q_1,q_2) \to \mathbb{R}$ denote the {\bfi stochastic HP action integral},
\[ 
\mathfrak{G}(\omega, q,v,p) = \int _a^b \left[ \mathcal{L}(q, v) dt + \sum_{i=1}^m \gamma_i (q) \circ d W_i(t,\omega) + \left\langle p, \frac{dq}{dt} - v \right\rangle dt \right] \text{.} 
\]   
\end{definition}

The Pontryagin bundle is a vector bundle over $Q$ whose fiber at $q \in Q$ is the vector space $T_q Q \oplus T_q^* Q$.   The HP path space is a smooth infinite-dimensional manifold.  One can show that its tangent space at  $c = (q,v,p) \in \mathcal{C}( [a,b], q_1, q_2)$ consists of maps $w = (q,v,p, \delta q, \delta v, \delta p) \in C^0([a,b],T(TQ \oplus T^*Q))$  such that  $\delta q(a) = \delta q(b) = 0$ and $q, \delta q$ are of class $C^1$.     
Let $(q,v,p)(\cdot, \epsilon)  \in \mathcal{C}(q_1, q_2, [a,b])$ denote a one-parameter family of curves in $\mathcal{C}$ that is differentiable with respect to $\epsilon$. Define the differential of $\mathfrak{G}$ as,
\[
\mathbf{d} \mathfrak{G} \cdot (\delta q, \delta v, \delta p) :=
\frac{\partial}{\partial \epsilon}  \left. \mathfrak{G}(\omega, q(t,\epsilon), v(t, \epsilon), p(t, \epsilon)) \right|_{\epsilon=0} 
\]
where 
\[
\delta q(t) = \frac{\partial}{\partial \epsilon} \left. q(t,\epsilon)  \right|_{\epsilon=0} ,~~\delta q(a) = \delta q(b) = 0,~~ 
\delta v(t) = \frac{\partial}{\partial \epsilon} \left. v(t,\epsilon)  \right|_{\epsilon=0} ,~~  
\delta p(t) = \frac{\partial}{\partial \epsilon} \left. p(t,\epsilon) \right|_{\epsilon=0} \text{.}
\]
In terms of this differential, one can state the following critical point condition.

\begin{theorem}[\textbf{Stochastic Variational Principle of Hamilton-Pontryagin}]
Let $\mathcal{L}: TQ \to \mathbb{R}$ be a Lagrangian on $TQ$ of class $C^2$ with respect to $q$ and $v$ and with globally Lipschitz first derivatives with respect to $q$ and $v$.  Let $\gamma_i: Q \to \mathbb{R}$ be of class $C^2$ and with globally Lipschitz first derivatives for $i=1,...,m$.  Then almost surely,
\begin{center}
\begin{quote}
A curve $c = (q,v,p) \in \mathcal{C}([a,b],q_1, q_2)$ satisfies the stochastic HP equations:
\begin{equation}
\label{eq:shp}
\begin{cases}
\begin{array}{rcl}
d q &=& v dt  \text{,}  \\
d p &= & \frac{\partial \mathcal{L}}{\partial q} dt + \sum_{i=1}^m \frac{\partial \gamma_i}{\partial q} \circ d W_i \text{,}  \\
p &=& \frac{\partial \mathcal{L}}{\partial v} \text{.} 
\end{array}
\end{cases}
\end{equation}
if and only if it is a critical point of the function $\mathfrak{G}: \Omega \times \mathcal{C}( [a,b], q_1, q_2) \to \mathbb{R}$, that is, $\mathbf{d} \mathfrak{G} (\omega, c(\omega)) = 0$. 
\end{quote}
\end{center}
\label{thm:stochastichpprinciple}
\end{theorem}

\begin{proof}
Let us first prove almost surely: a critical point of the function satisfies (\ref{eq:shp}).  The differential of the HP action integral is given by,
\begin{align*}
\mathbf{d} \mathfrak{G}  (\omega, c)  \cdot & (\delta q, \delta v, \delta p) = \int_a^b \left[ \frac{\partial \mathcal{L}}{\partial q} \cdot \delta q ds + \sum_{i=1}^m \frac{\partial \gamma_i}{\partial q} \cdot \delta q \circ  d W_i \right. \\
& \left. +  \frac{\partial \mathcal{L}}{\partial v} \cdot \delta v ds + \left\langle \delta p,  \frac{d q}{dt} - v  \right\rangle ds + \left\langle p, \delta \frac{d q}{dt} - \delta v  \right\rangle ds \right] \text{.}
\end{align*}
One can use a dominated convergence argument to show that differentiation and stochastic integration commute in the above stochastic integrals, as $\gamma_i$ are of class $C^2$  for $i=1,...,m$, and the curves are continuous.  Consider the term involving $\delta p$.  Since $\delta p$ is arbitrary and the integrand is continuous, the kinematic constraint holds: $dq/dt = v$.\footnote{This follows from the basic lemma that if $f, g \in C^0([a,b], \mathbb{R})$ and $g$ is arbitrary then $\int_a^b f(t) g(t) dt = 0 \iff f(t) = 0~~\forall~~t \in[a,b]$.}  Similarly, the Legendre transform is obtained from the $\delta v$ term: $\partial \mathcal{L}/\partial v = p$.

Collecting the variations with respect to $\delta q$ in the differential gives,
\[
\int_a^b \left[ \frac{\partial \mathcal{L}}{\partial q}  \cdot \delta q ds + \left\langle p, \delta \frac{d q}{dt} \right\rangle ds + 
  \sum_{i=1}^m \frac{\partial \gamma_i}{\partial q} \cdot \delta q \circ  d W_i   \right] \text{.}
\]
The first two terms are standard Lebesgue integrals and the last $m$ terms are Stratonovich stochastic integrals.   Motivated by Riemann-Stieljes integration by parts the following definition is introduced:

\begin{definition}
Given $f_1 \in C^0([a,b], T^*Q)$ and $f_2 \in C^1([a,b], TQ)$, define
\[
\int_a^t \left\langle df_1,  f_2 \right\rangle :=  
\left. \left\langle f_1, f_2 \right\rangle \right|_{a}^t - \int_a^t \left\langle f_1,  \frac{d f_2}{dt} \right\rangle dt,~~ t \in [a,b]   \text{.}
\]
\label{def:dfintegral}
\end{definition}

Using this definition and the boundary conditions $\delta q(a) = \delta q(b) = 0$,  the following function $I: \mathcal{C}(q_1, q_2, [a,b]) \times C^1([a,b], TQ)  \to \mathbb{R}$ is introduced,
\[
I(q,v,p,f) = \int_a^b \left[ \left( \frac{\partial \mathcal{L}}{\partial q} ds  + \sum_{i=1}^m \frac{\partial \gamma_i}{\partial q} \circ d W_i  - dp \right) \cdot f  \right]  \text{,}
\]
so that,
\[
I(q,v,p, \delta q) = \int_a^b \left[ \frac{\partial \mathcal{L}}{\partial q}  \cdot \delta q ds + \left\langle p, \delta \frac{d q}{dt} \right\rangle ds + 
  \sum_{i=1}^m \frac{\partial \gamma_i}{\partial q} \cdot \delta q \circ  d W_i   \right] \text{.}
\]
In the following it is shown that if $I(q,v,p, f) = 0$ for arbitrary $f$ of class $C^1$ then $(q,v,p)$ satisfy (\ref{eq:shp}).

Let $\{ U_{\alpha}, g_{\alpha} \}$ be a partition of unity on $TQ \oplus T^*Q$.  Expand $I$ in terms of this partition of unity,
\[
I = \sum_{\alpha} \int_a^b \left[ g_{\alpha}(q,v,p) \left( \frac{\partial \mathcal{L}}{\partial q} dt  + \sum_{i=1}^m \frac{\partial \gamma_i}{\partial q} \circ  d W_i  - dp \right) \cdot f \right]  \text{.}
\]
Since the curves $(q,v,p)$ are compactly supported, only a finite number of the $g_{\alpha}$ are nonzero.  For each $g_{\alpha}$ nonzero, the terms in the integral can be expressed in local coordinates.     Observe that since $dq = v dt$, the Stratonovish-Ito conversion formula implies that, 
\[
 \int_a^b  g_{\alpha} \frac{\partial \gamma_i}{\partial q} \cdot \delta q  \circ d W_i  =
  \int_a^b  g_{\alpha} \frac{\partial \gamma_i}{\partial q} \cdot \delta q  d W_i  
\]
for $i=1,...,m$.

We will select $f$ to single out the jth-component of the covector field in $I$.  Introduce the following function $h: \mathbb{R} \to \mathbb{R}$ for this purpose:
\[
h(t) =  2 \frac{t}{\epsilon} - \frac{t^2}{\epsilon^2} \text{.}
\] 
Observe that $h(0)=0$, $h(\epsilon) = 1$, and $h^{\prime}(\epsilon)=0$.  Let $\{ e_j \}$ be a basis for the model space of $Q$.  Now fix $j$, and define $f_{\epsilon} \in C^1([a,b], TQ)$ in local coordinates as follows:
\[
f_{\epsilon}(s)  = 
\begin{cases}
h(s-a) e_{j} &  \text{if}~ ~a \le s \le a + \epsilon \text{,} \\
e_{j}  &  \text{if}~~ a + \epsilon < s < t - \epsilon \text{,} \\
h(t-s) e_{j} &  \text{if}~~ t - \epsilon  \le s \le t  \text{,} \\
 0 &  \text{if}~~  t  < s \le b \text{.}
\end{cases}
\]
Introduce the following label to simplify subsequent calculations,
\[
A(s)  = \left( \frac{\partial \mathcal{L}}{\partial q}(q(s), v(s)) ds  + \sum_{i=1}^m \frac{\partial \gamma_i}{\partial q}(q(s)) ~d W_i(s) - dp(s) \right) \cdot e_j  \text{.}
\]
In terms of $A(s)$,  one can write
\[
I(q,v,p, f_{\epsilon}) = \sum_{\alpha} \left[ \int_{a}^{a+\epsilon} h(s-a) g_{\alpha}(s) A(s) 
+ \int_{a+\epsilon}^{t-\epsilon} g_{\alpha}(s) A(s) 
+ \int_{t-\epsilon}^t  h(t-s) g_{\alpha}(s) A(s) 
\right]  \text{.}
\]

The proof shows in the mean squared norm (cf.~definition~\ref{def:meansquared}),
\begin{equation}
\lim_{\epsilon \to 0}  I(q,v,p,f_{\epsilon})  = \sum_{\alpha} \int_{a}^{t} g_{\alpha} A(s) =: I^*  \text{.}
\label{eq:msIconvergence}
\end{equation}
Using this result and the Borel-Cantelli lemma, one can deduce there exists $\{ \epsilon_n \}$ that converges to $0$  such that $I(q,v,p,f_{\epsilon_n})$ a.s.~converges to $I^*$.  It follows that $I^* = 0$ almost surely.

We proceed to prove (\ref{eq:msIconvergence}).  Since $(a+b)^2 \le 2 a^2 + 2 b^2$,
\begin{align*}
 & \left\|  \sum_{\alpha} \int_{a}^{t} g_{\alpha} A(s) -  I(q,v,p,f_{\epsilon})  \right\|^2  \\
& = \left\|  \sum_{\alpha}
\int_{a}^{a+\epsilon} ( 1- h(s-a)) g_{\alpha} A(s) 
+ \int_{t-\epsilon}^t (1- h(t-s)) g_{\alpha} A(s) \right\|^2  \\
& \le   2 \left\| \sum_{\alpha} \int_{a}^{a+\epsilon} (1-h(s-a)) g_{\alpha} A(s) \right\|^2 
+ 2 \left\| \sum_{\alpha} \int_{t-\epsilon}^t  (1-h(t-s))  g_{\alpha} A(s) \right\|^2
\text{.}
\end{align*}
We will only show how to bound the first term since bounding the second term is very similar.  By continuity of $(q,v,p)$, one can pick $\epsilon$ small enough so that the support of $(q,v,p)$ lies in a single chart.  Therefore,
\begin{align*}
 & \left\|  \sum_{\alpha} \int_{a}^{a+\epsilon} (1-h(s-a)) g_{\alpha} A(s) \right\|^2  \\
& = \left\| \int_{a}^{a+\epsilon} (1-h(s-a))  
\left( \frac{\partial \mathcal{L}}{\partial q} ds  + \sum_{i=1}^m \frac{\partial \gamma_i}{\partial q} ~d W_i - dp \right) \cdot e_j  \right\|^2  \\
&\le 3 \left\| \int_{a}^{a+\epsilon}  (1-h(s-a))  \frac{\partial \mathcal{L}}{\partial q^j} ds \right\|^2   +  3 \sum_{i=1}^m \left\| \int_{a}^{a+\epsilon} (1-h(s-a)) \frac{\partial \gamma_i}{\partial q^j} ~d W_i \right\|^2  \\
 &+ 3 \left\|   \int_{a}^{a+\epsilon} (1-h(s-a)) dp \cdot e_j \right\|^2 \text{.}
\end{align*}
Since $\frac{\partial \mathcal{L}}{\partial q^j}$ is continuous on $s \in[a,a+\epsilon]$, the first term can be bounded,
\[
\left\| \int_{a}^{a+\epsilon} (1-h(s-a))   \frac{\partial \mathcal{L}}{\partial q^j} ds \right\|^2 \le
\frac{M^2  \epsilon^2}{9} \text{.}
\]
Similarly, by the Ito isometry and since $\frac{\partial \gamma_i}{\partial q^j}$ is continuous on $s \in[a, a +\epsilon]$, the second $m$ terms can similarly be bounded, e.g., the ith Stratonovich integral can be bounded as follows,
\[
 \left\| \int_{a}^{a+\epsilon} (1- h(s-a)) \frac{\partial \gamma_i}{\partial q^j} ~d W_i \right\|^2  =
  E \left( \int_{a}^{a+\epsilon} \left| (1-h(s-a)) \frac{\partial \gamma_i}{\partial q^j} \right|^2 ds \right) \le \frac{M^2 \epsilon}{5} \text{.}
\]
Using definition~\ref{def:dfintegral} and the integral mean value theorem, the final term can be bounded as well:
\begin{align*}
\left\|   \int_{a}^{a+\epsilon} (1-h(s-a)) dp \cdot e_j \right\|^2 & = 
 \left\| \left. (1-h(s-a)) p_j (s) \right|_{s=a}^{s=a+\epsilon} + \int_a^{a+\epsilon} p_j(s) h^{\prime}(s-a) ds \right\|^2  \\
 & =   \left\| - p_j(a)  + p_j(a + c \epsilon)  \right\|^2  
\end{align*}
where $0 \le c \le 1$ is a real constant.  Since $p_j$ is of class $C^0$, as $\epsilon \to 0$ this term vanishes.  Since $j$ is arbitrary we have proved (\ref{eq:msIconvergence}).  Therefore, almost surely: if $c$ is a critical point of $\mathfrak{G}$ then $\mathbf{d}  \mathfrak{G} (\omega, c)  \cdot w =0 $ for all $w \in T_c  \mathcal{C}( [a,b], q_1, q_2) $, and hence, $c$ satisfies the stochastic HP equations.

On the other hand, almost surely: if $c$ satisfies  (\ref{eq:shp}), then it is a critical point of $\mathfrak{G}$.  This direction is easy to confirm, since as a solution to the stochastic HP equations $c$ is a measureable diffusion process.   In fact, this direction is similar to the one Bismut originally established, namely that the solution of stochastic Hamilton's equations extremize an action function; albeit the stochastic action used by Bismut has a different domain than the stochastic action used in this proof \citep{Bi1981}.  
\end{proof}

Equations (\ref{eq:shp}) are a stochastic differential algebraic system of equations. Assuming one can eliminate $v$ using the Legendre transform, these equations can be viewed as a Cauchy problem. This paper is primarily concerned with forces or torques that appear as white noise in the balance of momentum equations, which explains the choice $\gamma_i=\gamma_i(q)$.  Observe that by the Ito-Stratonovich conversion formula, the Ito modification to the drift is equal to zero, and hence, (\ref{eq:shp}) can be written in Ito form as:
\begin{align*}
d q &= v dt  \text{,} \\
d p &=  \frac{\partial \mathcal{L}}{\partial q} dt + \sum_{i=1}^m \frac{\partial \gamma_i}{\partial q} d W_i \text{,}  \\
p &= \frac{\partial \mathcal{L}}{\partial v} \text{.} 
\end{align*}
In what follows structure-preserving properties of the flow map defined by the maximal solution of these equations over $[a,b]$ will be investigated.   First, observe that because of the smoothness conditions assumed in theorem~\ref{thm:stochastichpprinciple}, a solution almost surely exists and is pathwise unique on $[a,b]$ by the results in \S 2.  When $\gamma_i$ is constant for $i=1,...,m$,  the reader is referred to the following texts for deterministic treatments  of symplecticity, momentum map preservation, and holonomically constrained mechanical systems: \citep{MaRa1999, MaWe2001}.

 \paragraph{Symplecticity}

Assuming that one can eliminate $v$ using the Legendre transform, the stochastic HP equations define a {\em stochastic flow} on the symplectic manifold $(T^*Q, \kappa)$ where $\kappa$ is the canonical symplectic form \citep{MaRa1999}.   We denote this flow by $F_t: T^*Q \to T^*Q$.   With this assumption and in a more general context, Bismut extended the variational proof of symplecticity and Noether's theorem to stochastic Hamiltonian systems \citep{Bi1981}.   In fact, one does not need to prove both directions of (\ref{thm:stochastichpprinciple}) to perform this extension.  These proofs are repeated here in the context of stochastic Hamiltonian systems driven by Wiener processes for the reader's convenience and for completeness.

The variational proof of symplecticity will be used to show this flow preserves $\kappa$ \citep{MaRa1999}.  By theorem (\ref{thm:maximalsolution}), assuming the Lagrangian is sufficiently smooth  then $F_t: T^*Q \to T^*Q$ is a.s.~differentiable.  Recall the derivative is defined in the mean squared sense (cf.~definition~\ref{def:meansquared}).     The idea of the proof is to restrict $\mathfrak{G}$ to the space of pathwise unique solutions, i.e., define $\hat{\mathfrak{G}} = \left. \mathfrak{G} \right|_{\text{solutions}}$.  On the same filtered probability space with the same Brownian motion, this solution space can be identified with the set of initial conditions, i.e., this restricted action can be expressed as, $\hat{\mathfrak{G}}(\omega, \cdot) : T^*Q \to \mathbb{R}$.  For each initial condition by  (\ref{thm:maximalsolution}), there exists a pathwise unique solution almost surely.  One then computes $\mathbf{d} \hat{\mathfrak{G}}$:
\begin{align*}
\mathbf{d} & \hat{\mathfrak{G}} (\omega, q(a), p(a))  \cdot (\delta q(a), \delta p(a)) = \\
&\int_a^b \left[ \left(  \frac{\partial \mathcal{L}}{\partial q} dt + \sum_{i=1}^m \frac{\partial \gamma_i}{\partial q} \circ d W_i  - dp \right) \cdot \delta q  + \delta p \cdot \left( \frac{dq}{dt} - v  \right) ds + 
\left( \frac{\partial \mathcal{L}}{\partial v}  - p \right) \cdot \delta v dt \right] \\
&+ \left. \left\langle p, \delta q \right\rangle \right|_{a}^b   \text{.}
\end{align*}
The integral in the above vanishes since $\hat{\mathfrak{G}}$ is restricted to solution space.  The boundary terms define in local coordinates the canonical 1-form $\Theta$ on $T^*Q$.  Computing $\mathbf{d}^2 \hat{\mathfrak{G}}$ gives conservation of the canonical symplectic form.

\begin{theorem}[\textbf{Conservation of stochastic symplectic form}]
Provided that one can eliminate $v$ using the Legendre transform, the flow of (\ref{eq:shp}) preserves the canonical symplectic form almost surely.  
\end{theorem}

 \paragraph{Noether's Theorem}

In what follows we review for completeness Bismut's extension of Noether's theorem \citep{Bi1981}.   Let $G$ be a Lie group with Lie algebra $\mathfrak{g}$.   The {\it left action} of the Lie group on  $Q$ is denoted $\Phi: G \times Q \to Q$.  The {\it cotangent lift}  of this action is likewise denoted $ \Phi^{T^*Q}: G \times T^*Q \to T^*Q  $: 
\begin{align*}  
\Phi^{T^*Q}(h, q,p)  &= 
\left( \Phi(h, q), ((D_q \Phi(h, q))^{-1})^* \cdot p \right)  \text{.}
\end{align*}   
The corresponding {\it infinitesimal generators} are $\xi^Q: Q \to T^*Q$ and $\xi^{ T^*Q }:  T^*Q \to T( T^*Q)$ and by definition ,
\[
\xi^Q(q,p) = \frac{d}{ds} \left[ \Phi(\exp(s \xi), q) \right]_{s=0} \text{,}~~~  
\xi^{T^*Q}(q, p) = 
\frac{d}{ds} \left[ \Phi^{ T^*Q }(\exp(s \xi), q, p) \right]_{s=0}   \text{.}
\]

This action gives rise to the following momentum map $J: T^*Q  \to \mathfrak{g}^*$  
\[
J(q, p) \cdot \xi = \left\langle p,  \xi^Q(g, p)  \right\rangle \text{,}
\]  
The following conservation law follows if $\mathfrak{G}$ is infinitesimally invariant with respect to the G-action.   We remark in passing that infinitesimal invariance of $\mathfrak{G}$ follows from left-invariance of the stochastic HP action with respect to the G-action.

\medskip

\begin{theorem}[\textbf{Stochastic Noether's Theorem}] Let $G$ be a Lie group.  If $\mathfrak{G}$ is infinitesimally symmetric with respect to the left (or right) action of $G$, then the corresponding momentum map is conserved a.s., i.e., $J = \left\langle p, \xi_Q(q) \right\rangle$, is a conserved quantity under the flow of (\ref{eq:shp}).
\end{theorem}

\begin{proof}  
Consider the differential of $\hat{\mathfrak{G}}$ in the direction of $\xi^{T^*Q}$,
\begin{align*}
\mathbf{d} \hat{\mathfrak{G}} (\omega, q(a), p(a))  \cdot \xi^{T^*Q}(q(a), p(a)) =  \left. \left\langle p, \xi^Q(q,p) \right\rangle \right|_{a}^b   \text{.}
\end{align*}
Moreover infinitesimal symmetry implies that
 \begin{align*}  
 \mathbf{d} \hat{\mathfrak{G}}  \cdot \xi^{T^*Q}(q(a), p(a))  &= 0 \text{,} \\
\implies J (q(b), p(b)) \cdot \xi  - J(q(a), p(a)) \cdot \xi &= 0  \text{,}
\end{align*}
and hence $J$ is conserved under the flow since $\xi$ is arbitrary.  
\end{proof}

\paragraph{Holonomic Constraints}

The following results will require the proof of the converse of theorem~\ref{thm:stochastichpprinciple}.  The setting in this part is an n-manifold $Q$ and a stochastic Hamiltonian system with holonomic constraint.   To be specific, suppose that the motion of the mechanical system is constrained to a submanifold $S \subset Q$ defined as $S = g^{-1}(0)$ where $g: Q \to \mathbb{R}^k$, $k < n$, $g$ is smooth, and $0$ is a regular value of $g$.

As opposed to using generalized coordinates on $TS$, we wish to describe the mechanical system using constrained coordinates on $TQ$ and introduce Lagrange multipliers to enforce the constraint.  However, because of the stochastic component of the action, the standard Lagrange multiplier theorem will not apply directly and one cannot introduce Lagrange multipliers in the standard way.  Instead, we will introduce the Lagrange multiplier using definition~\ref{def:dfintegral}.

In particular, consider the following constrained variational principle,
\[
\delta \left( \mathfrak{G} + \int_a^b \left\langle d \lambda, g \right\rangle \right) = 0 \text{,}
\]  
where using definition~\ref{def:dfintegral} 
\[
\int_a^t \left\langle d \lambda,  g \right\rangle :=  
\left. \left\langle \lambda, g \right\rangle \right|_{a}^t - \int_a^t \left\langle \lambda,  \frac{d g}{dt} \right\rangle dt,~~ t \in [a,b]   \text{.}
\]
In this case $\lambda(t)$ is a Lagrange multiplier dual to $g(q(t))$ for $t\in[a,b]$, and we assume is of class $C^0$.    The corresponding equations of motion are obtained in a similar fashion as (\ref{eq:shp}).  In particular, one can prove the following equivalence:

\begin{theorem}[{\bf Constrained, Stochastic HP Principle}]  Given a stochastic Hamiltonian system with Lagrangian $\mathcal{L}: TQ \to \mathbb{R}$ such that $\partial^2 \mathcal{L}/ \partial v^2$ is invertible, stochastic potentials $\gamma_j: Q \to \mathbb{R}$ for $j=1,...,m$, and holonomic constraint $g: Q \to \mathbb{R}^k$ with $S = g^{-1}(0)$. Set $\mathcal{L}^S = \left.\mathcal{L} \right|_{TS}$ and $\gamma_j^S = \left. \gamma_j \right|_{S}$ for $j=1,...,m$.   Fix $q_1, q_2 \in S$ and let $i: S \to Q$ be the inclusion mapping.  Then the following are equivalent:
\begin{description}
\item[(i)] Stochastic HP principle for $\mathcal{L}^S$ and $\gamma_j^S$, $j=1,...,m$, on $TS \oplus T^*S$ (cf.~thm.~\ref{thm:stochastichpprinciple}).
\item[(ii)] There exists $\lambda \in C^0([a,b],\mathbb{R}^k)$ such that $z=(q, v, p) \in \mathcal{C}([a,b],i(q_1),i(q_2))$ and $\lambda$ extremize the augmented action $\bar{\mathfrak{G}}(z,\lambda) =  \mathfrak{G}(z) + \left\langle \lambda, \Phi(z) \right\rangle$  where $\Phi(z)(t) = dg/dt(q(t))$ and $\left\langle \lambda, \Phi(z) \right\rangle = \int_a^b \left\langle \lambda(s), \Phi(z)(s) \right\rangle ds$.   
\item[(iii)] There exists $\lambda \in C^0([a,b],\mathbb{R}^k)$ such that $z=(q, v, p) \in \mathcal{C}([a,b],i(q_1),i(q_2))$ and $\lambda$ satisfy the constrained, stochastic HP equations
\begin{equation} \label{eq:cstochastichp}
\begin{cases}
\begin{array}{ccl}
dq &=& v dt \text{,}  \\
d p &=&  \frac{\partial \mathcal{L}}{\partial q}(q,v) dt + \sum_{j=1}^m \frac{\partial \gamma_j}{\partial q}(q) \circ d W_j + \frac{\partial g}{\partial q}(q)^* \cdot d \lambda \text{,}  \\
p &=& \frac{\partial \mathcal{L}}{\partial v}(q,v) \text{,}   \\
g(q) &=& 0  \text{.}  
\end{array}
\end{cases}
\end{equation}
\end{description}
\label{thm:stochasticconstrainedhpprinciple}
\end{theorem}

From this equivalence it follows that the flow of (\ref{eq:cstochastichp}) is mean-square symplectic.  For a proof of this theorem the reader is referred to \citep{BoOw2007b}.

\paragraph{Nonconservative Effects}

Nonconservative effects are incorporated by considering the {\bfi Lagrange-d'Alembert-Pontryagin principle}.  In this principle the effect of a nonconservative force appears as virtual work.  Consider a force field $F: TQ \to T^*Q$.  Then the Lagrange-d'Alembert-Pontryagin principle is given by
\[
\delta \int _a^b \left[\mathcal{L}(q, v) dt + \sum_{i=1}^m \gamma_i (q) \circ d W_i + \left\langle p, \frac{dq}{dt} - v \right\rangle dt \right]  + \int_a^b F(q, v) \cdot \delta q dt  = 0,
\]  
where the variations of the first term are the usual ones (vanishing at the endpoints).  This principle provides a simple way to add dissipative effects into the drift which, e.g., appear in the standard Langevin equations for particles.

\paragraph{Lagrangian Reduction}  For background and exposition on Lagrangian reduction in the deterministic setting the reader is referred to \citep{MaSc1993, MaRa1999}.  Suppose that $Q$ is a  Lie group $G$ with Lie algebra $\mathfrak{g}$.  In this context one can define a {\em left-trivialized Lagrangian} by using the left-action of $G$ to left-trivialize $\mathcal{L}$.  One does this by transforming a point $(g(t), v(t)) \in TG$ to $(g(t), \xi(t)) \in G \times \mathfrak{g}$ via the relation between the velocity at $g(t) \in G$ and the {\em body angular velocity} at $e \in G$ given by: $\xi(t)= g(t)^{-1} v(t) \in \mathfrak{g}$.   Denote by $l:G \times \mathfrak{g} \to \mathbb{R}$ the deterministic left-trivialized Lagrangian defined as $l(g(t), \xi(t)) = \mathcal{L}( g(t), g(t) \xi(t))$. The variational principle associated with $l$ is the {\bfi left-trivialized HP principle} which can be written as:
\[
\delta \int_a^b \left[ l(g,\xi) dt + \sum_{i=1}^m \gamma_i(g) \circ d W_i  + \left\langle \mu, g^{-1} \frac{dg}{dt} - \xi  \right\rangle dt  \right] = 0
\]
In this principle the Lagrange multiplier $\mu \in \mathfrak{g}^*$ is the body angular momentum. For more details on the geometry of this principle in the deterministic setting the reader is referred to \citep{BoMa2007, Bo2007}.   The resulting equations are obtained by taking arbitrary variations with fixed endpoints on $g$. For a function $U: G \to \mathbb{R}$, define its {\bfi left-trivialized differential}, $U_g \in \mathfrak{g}^*$, as
\begin{equation}
U_g \cdot \eta := \left\langle  \frac{\partial U}{\partial g}, TL_g \eta \right\rangle, ~~~ \eta \in \mathfrak{g} \text{.}
\label{def:lefttrivializeddifferential}
\end{equation}
In terms of $U_g$, one can write the {\bfi stochastic left-trivialized HP equations}:
\begin{align}
\frac{d g}{dt} &= g \xi   \label{eq:rhp1}, \\
d \mu &= l_g dt + \sum_{i=1}^m (\gamma_i)_g \circ d W_i(\omega, t) \label{eq:rhp2}, \\
\mu &= \frac{\partial l}{\partial \xi} \label{eq:rhp3} \text{.}
\end{align}
By eliminating $\xi$ using (\ref{eq:rhp3}), one obtains a SDE on $G \times \mathfrak{g}^*$.  The kinematic constraint in this context is referred to as the {\em reconstruction equation}.  We summarize this section with the following theorem.

\begin{theorem}[\textbf{Stochastic Left-Trivialized Hamilton-Pontryagin Principle}]  Consider a mechanical system on a Lie group $G$ with left-trivialized Lagrangian $l: G \times \mathfrak{g} \to \mathbb{R}$ and stochastic potentials $\gamma_i: G \to \mathbb{R}$ for $i=1,...,m$.  Let $s$ denote the left-trivialzied action given by,
\[
s = \int_a^b \left[  l(g,\xi) dt + \sum_{i=1}^m \gamma_i(g) \circ d W_i + \left\langle \mu, g^{-1} \frac{dg}{dt} - \xi  \right\rangle dt  \right]  \text{.}
\]
Almost surely, the stochastic HP principle on a Lie group (cf.~thm.~\ref{thm:stochastichpprinciple}), is equivalent to the stochastic left-trivialized HP principle:
\[
\delta s = 0,
\]
where the curves
\[
g(t) \in G, \xi(t) \in \mathfrak{g}, \mu(t) \in \mathfrak{g}^*,~~t\in[a,b]
\]
can be varied arbitrarily with $\delta g(a) = \delta g(b) = 0$.  A curve is a critical point of the left-trivialized action if and only if it satisfies the left-trivialized HP equations given by (\ref{eq:rhp1})-(\ref{eq:rhp3}).  
\end{theorem}

\section{Stochastic Variational Integrators}  \label{sec:stochasticvi}

The standard approach of deriving variational integrators is extended to the stochastic context in this section; see, e.g., \citep{MaWe2001}.  The cornerstone of variational integration theory is the discrete Lagrangian.    In this section we develop and analyze integrators from an abstract discrete Lagrangian which takes values on the configuration space squared.  In the subsequent sections, discrete Lagrangians will be specified and the associated time-integrators analyzed from the HP viewpoint.    Let $[a,b]$ and $N$ be given, and define the fixed step size $h=(b-a)/N$ and $t_k = h k + a$, $k=0,...,N$.

\begin{definition}
Consider a mechanical system with given discrete Lagrangian $\mathcal{L}_d: Q \times Q \to \mathbb{R}$.   Let  $\theta_t: \Omega \to \Omega$, $t \in [a,b]$ denote the base flow on the probability space (cf.~definition~\ref{def:stochasticflow}). Let the approximant to the Stratonovich integrals be denoted by:
\[
B_d(t_k, q_k, q_{k+1}, \omega) h \approx \sum_{i=1}^m \int_{t_k}^{t_{k+1}} \gamma_i(q(t)) \circ d W_i(t, \theta_{t_k+t}(\omega))
\]
The associated {\bfi stochastic discrete Lagrangian} $L_d: \mathbb{R} \times \Omega \times Q \times Q \to \mathbb{R}$ is defined as:
\[
L_d(t_k, \omega, q_k, q_{k+1}) = \mathcal{L}(q_k, q_{k+1}) +  B_d(t_k, q_k, q_{k+1}, \omega) \text{.}
\]
Fixing the interval $[a,b]$, define the {\bfi discrete path space}  as
\[
\mathcal{C}_d(Q) = \{ q_d: \{ t_k \}_{k=0}^N \to Q \}  \text{.}
\] 
Let $\mathfrak{G}_d: \Omega \times \mathcal{C}_d(Q) \to \mathbb{R}$ denote the {\bfi stochastic action sum},
\[ 
\mathfrak{G}_d(\omega, q_d) = \sum_{k=0}^{N-1} L_d(t_k, \omega, q_k, q_{k+1}) h  \text{.}
\]   
\end{definition}
The {\bfi discrete stochastic Hamilton's principle} states that the path the mechanical system takes in $\mathcal{C}_d$ is one that extremizes $\mathfrak{G}_d(\omega, \cdot)$ subject to fixed endpoint conditions $q_0=q(a)$ and $q_N=q(b)$.  By discrete integration by parts (re-indexing),
\begin{align*}
\mathbf{d} \mathfrak{G}_d(\omega, q_d) \cdot \{ \delta q_k \} =& 
\sum_{k=1}^{N-1}
\left( D_3 L_d (t_k, \omega, q_k, q_{k+1}) + D_4 L_d( t_{k-1}, \omega, q_{k-1}, q_{k}) \right) \cdot \delta q_k \\
&+ D_3 L_d (t_0, \omega, q_0, q_1) \cdot \delta q_0 + D_4 L_d (t_{N-1}, \omega, q_{N-1}, q_N) \cdot \delta q_N \text{.}
\end{align*}
Using the endpoint conditions $\delta q_0 = \delta q_N = 0$ one obtains:
\[
\mathbf{d} \mathfrak{G}_d(\omega, q_d) \cdot \{ \delta q_k \} = 
\sum_{k=1}^{N-1}
\left( D_3 L_d (t_k, \omega, q_k, q_{k+1}) + D_4 L_d( t_{k-1}, \omega, q_{k-1}, q_{k}) \right) \cdot \delta q_k \text{.}
\]
Stationarity of this action sum implies the following {\bfi stochastic discrete Euler-Lagrange} equations:
\[
D_3 L_d (t_k, \omega, q_k, q_{k+1}) + D_4 L_d( t_{k-1}, \omega, q_{k-1}, q_{k}) = 0 
\]
for $k=1, ... , N-1$.  The resulting update scheme is not self-starting.  To initialize the method one needs to provide $(q_0, q_1) \in Q \times Q$ as opposed to a point $(q_0, v_0) \in TQ$.

\paragraph{Symplecticity}
As in the continuous theory symplecticity follows from restricting $\mathfrak{G}_d(\omega, \cdot)$ to pathwise unique solutions of the stochastic discrete Euler-Lagrange equations, $\hat{\mathfrak{G}}_d$.  Since pathwise unique solutions can be parametrized by initial conditions we regard the restricted action as $\hat{\mathfrak{G}}_d: \Omega \times Q \times Q \to \mathbb{R}$.  Taking its first variation gives,
\begin{align*}
\mathbf{d} \hat{\mathfrak{G}}_d(\omega, q_0, q_1) &\cdot ( \delta q_0, \delta q_1 ) = 
\sum_{k=1}^{N-1}
\left( D_3 L_d (t_k, \omega, q_k, q_{k+1}) + D_4 L_d( t_{k-1}, \omega, q_{k-1}, q_{k}) \right) \cdot \delta q_k \\
&+ D_3 L_d (t_0, \omega, q_0, q_1) \cdot \delta q_0 + D_4 L_d (t_{N-1}, \omega, q_{N-1}, q_N) \cdot \delta q_N \text{.}
\end{align*}
Because of the restriction to solution space the sum vanishes, and the boundary terms remain:
\[
\mathbf{d} \hat{\mathfrak{G}}_d(\omega, q_0, q_1) \cdot ( \delta q_0, \delta q_1 ) =
D_3 L_d (t_0, \omega, q_0, q_1) \cdot \delta q_0 + D_4 L_d (t_{N-1}, \omega, q_{N-1}, q_N) \cdot \delta q_N \text{.}
\]
These boundary terms define left and right one-forms as follows:
\begin{align*}
\Theta^+(t_k, \omega, q_k, q_{k+1})  \cdot (\delta q_k, \delta q_{k+1})  &= D_4 L_d (t_k, \omega, q_k, q_{k+1}) \cdot \delta q_{k+1} \\
\Theta^-(t_k, \omega, q_k, q_{k+1}) \cdot (\delta q_k, \delta q_{k+1}) &= D_3 L_d (t_k, \omega, q_k, q_{k+1}) \cdot \delta q_k 
\end{align*}
which from $\mathbf{d}^2 L_d = 0$ satisfy:
\[
\mathbf{d} \Theta^+ = - \mathbf{d} \Theta^- = \Omega \text{.}
\]
Applying the second exterior derivative to $\hat{\mathfrak{G}}_d$ implies,
\[
\Omega(t_{N-1}, \omega, q_{N-1},q_N)(\delta q_{N-1}^1, \delta q_N^1) (\delta q_{N-1}^2, \delta q_N^2) = 
\Omega(t_0, \omega, q_{0},q_1)(\delta q_0^1, \delta q_1^1) (\delta q_0^2, \delta q_1^2) 
\]
since $\mathbf{d}^2 \hat{\mathfrak{G}}_d = 0$.  Hence, the discrete flow preserves the symplectic form $\Omega$.

\paragraph{Discrete Momentum Map}  Consider the left action of a Lie Group $G$ on $Q$.  If the stochastic discrete Lagrangian is infinitesimally symmetric, then the associated momentum map is preserved.  A sufficient condition for this is that the discrete Lagrangian is invariant with respect to the left action of $G$.  The proof is sketched out here.

Let the action on the discrete configuration manifold be denoted by $\Phi_{Q \times Q}: G \times Q \times Q \to Q \times Q$ and defined by:
\[
\Phi_{Q \times Q}(g, q_1, q_2) = (\Phi(g, q_1), \Phi(g, q_2)) \text{.}
\]
The associated infinitesimal generator is denoted $\xi_{Q \times Q}: Q \times Q \to T(Q \times Q)$ and by definition
\[
\xi_{Q \times Q}(q_1, q_2) = \frac{d}{ds} \left. \Phi_{Q \times Q}(\exp(s \xi), q_1, q_2) \right|_{s=0} \text{.}
\]
Assume that $L_d$ is infinitesimally symmetric, i.e., 
\[
\mathbf{d} L_d \cdot \xi_{Q \times Q} = \Theta^+ \cdot \xi_{Q \times Q} + \Theta^- \cdot \xi_{Q \times Q} = 0 \text{.}
\]
By this condition the left and right discrete momentum maps, $J^+, J^-: Q \times Q \to \mathfrak{g}^*$,:
\begin{align*}
J^+ \cdot \xi &= \Theta^+ \cdot \xi_{Q \times Q} \\
J^- \cdot \xi  &= - \Theta^- \cdot \xi_{Q \times Q}
\end{align*}
are equal, i.e., $J^+ = J^- = J$.  Consider the restricted action sum and compute its differential in the direction of the infinitesimal generator to obtain:
\begin{align*}
\mathbf{d} \hat{\mathfrak{G}}_d(\omega, q_0, q_1) & \cdot \xi_{Q \times Q}( q_0, q_1 ) =
\Theta^-(t_0, \omega, q_0, q_1) \cdot \xi_{Q \times Q}( q_0,q_1 )  \\
&+ \Theta^+(t_{N-1}, \omega, q_{N-1}, q_N) \cdot \xi_{Q \times Q}( q_{N-1}, q_N ) \text{,}
\end{align*}
which can be rewritten in terms of the momentum maps as
\begin{align*}
\mathbf{d} \hat{\mathfrak{G}}_d(\omega, q_0, q_1)  \cdot \xi_{Q \times Q}( q_0, q_1 ) =
-J^-(t_0, \omega, q_0, q_1) \cdot \xi 
+ J^+(t_{N-1}, \omega, q_{N-1}, q_N) \cdot \xi = 0 \text{.}
\end{align*}
Since the left and right momentum maps evaluated at the same point are equal, the momentum map $J$ is preserved under the discrete flow.

\section{Stochastic Variational Euler Integrator} \label{sec:stochasticveuler}

In the deterministic setting, the HP context provides a practical way to design discrete Lagrangians and obtain one-step methods on $TQ$ or $T^*Q$ as pointed out in \citep{Bo2007}.  In this section we examine variational Euler methods extended to the stochastic context following the  continuous stochastic HP theory laid out in \S 2.

\paragraph{Stochastic Variational Euler on $\mathbb{R}^n$}

To discretize the stochastic HP action integral one needs to replace the continuous Lagrangian, stochastic integral and kinematic contraint by discrete approximants.   We begin by introducing a first-order discretization of the kinematic constraint in (\ref{eq:shp}).  Let $[a,b]$ and $N$ be given and define the fixed step size $h=(b-a)/N$ and $t_k = h k$, $k=0,...,N$.

A discretization of the kinematic constraint can be obtained by introducing a discrete sequence $\{ q_k \}_{k=0}^N$ such that $q_k \in Q$ and  a finite difference map $ \varphi: Q \times Q \to TQ$.  For example, if $Q$ is a vector space the following backward difference map can be introduced:
\[
\varphi(q_k, q_{k+1}) = \left(q_{k+1}, \frac{q_{k+1}-q_k}{h}  \right) \text{.}
\]
Let $B_i^k \sim \mathcal{N}(0,h)$ be normally distributed random variables for $i=1,...,m$ and $k=0,...,N-1$.  In terms of the discretization of the kinematic constraint, the corresponding discrete HP action sum takes the following simple form:
\[
 \mathfrak{G}_d^e = \sum_{k=0}^{N-1}  \left[ \mathcal{L}(q_k, v_k) h +  
 \sum_{i=1}^m \gamma_i(q_k) B_i^k +
 \left\langle p_{k+1}, (q_{k+1}-q_k)/h - v_{k+1} \right\rangle h \right] 
\]
The stochastic discrete HP equations are given by:
\begin{align*}
q_{k+1} &=  q_k + h v_{k+1} \text{,} \\
p_{k+1} &= p_k + h \frac{\partial \mathcal{L}}{\partial q}(q_k, v_k) +  \sum_{i=1}^m \frac{\partial \gamma_i}{\partial q}(q_k) B_i^k \text{,} \\ 
p_k &= \frac{\partial \mathcal{L}}{\partial v}(q_k, v_k) \text{.}
\end{align*}

\paragraph{Stochastic Variational Euler on Lie Groups}

In the context of Lie groups, the reconstruction equation is discretized using canonical coordinates of the first kind, $\tau : \mathfrak{g} \to G$, as explained in \citep{Bo2007, BoMa2007}.   As in the vector space case, we define a finite difference map $ \varphi: G \times G \to G \times \mathfrak{g}$ that provides a first-order approximant to the reconstruction equation:
\[
\varphi(g_k, g_{k+1}) = (g_{k+1}, \tau^{-1} (g_k^{-1} g_{k+1})/ h) \in G \times \mathfrak{g} \text{.}
\]
A first order approximant to the stochastic left-trivialized action integral is given by, 
\begin{equation}
s_d=  \sum_{k=0}^{N-1}  \left[ l(g_k, \xi_{k}) h +  
 \sum_{i=1}^m \gamma_i(g_k) B_i^k  + 
 \left\langle \mu_{k+1}, \tau^{-1}(g_k^{-1} g_{k+1})/h - \xi_{k+1} \right\rangle h \right]   \text{.} \label{eq:reducedhpsum}
 \end{equation} 
 Let $d \tau^{-1}: \mathfrak{g} \times \mathfrak{g} \to \mathfrak{g}$ denote the right trivialized tangent of $\tau^{-1}$ as defined in \citep{BoMa2007}. The stochastic left-trivialized discrete HP equations are (cf.~definition~\ref{def:lefttrivializeddifferential}):
\begin{align*}
g_{k+1} &=  g_k \tau(h \xi_{k+1} )  \text{,} \\
(d \tau^{-1}_{h \xi_{k+1}})^* \mu_{k+1} &= (d \tau^{-1}_{-h \xi_k})^* \mu_k + h l_g(g_k, \xi_k) +  \sum_{i=1}^m (\gamma_i)_g(g_k) B_i^k \text{,} \\ 
\mu_k &= \frac{\partial l}{\partial \xi}(g_k, \xi_k)  \text{.}
\end{align*}

\paragraph{Holonomic Constraints \& Nonconservative Effects}

Holonomic constraints can be added via discrete Lagrange multipliers, and nonconservative effects via discrete impulses as described below in the Lie group context.    Suppose that $G$ is an n-manifold and that the mechanical system evolves on a submanifold $S \subset G$ defined as the zero level-set of $\varphi: G \to \mathbb{R}^k$ where $k < n$ and $S = \varphi^{-1}(0)$.  Further suppose that there exists a force field $F: G \to T^* G$.  These effects are appended by using the following modified stochastic left-trivialized HP action principle:
\begin{align*}
  \delta  \sum_{k=0}^{N-1}  & \left[ l(g_k, \xi_{k}) h  +  
 \sum_{i=1}^m \gamma_i(g_k) B_i^k  + 
 \left\langle \mu_{k+1}, \tau^{-1}(g_k^{-1} g_{k+1})/h - \xi_{k+1} \right\rangle h + 
  \left\langle \lambda_k, \varphi(g_k) \right\rangle h \right]  \\
 & +  \sum_{k=0}^{N-1}  F(g_k) \cdot \delta g_k h    = 0 \text{.} 
 \end{align*} 
 This algorithm for $G=\mathbb{R}^n$ is the stochastic analog of constrained symplectic Euler, and the numerical analysis of this method is discussed in \citep{BoOw2007b}.

\section{Langevin-type equations for Multiple Rigid Bodies}  \label{sec:rigidbodieslangevin}

\paragraph{Continuous Description}

Consider a mechanical system consisting of $K$ rigid bodies interacting via a pairwise potential dependent on their positions and orientations.  Let $(x_i (t), v_i (t), R_i (t), \omega_i (t)) \in T SE(3)$ denote the translational position, translational velocity, orientation, and spatial angular velocity of body $i$ where $i=1,...,K$.  Let $m_i$ and $\mathbb{I}_i$ denote the mass and diagonal inertia tensor of body $i$.  The left-trivialized Lagrangian for the system is given by:
\[
l(x_i, v_i, R_i, \omega_i) = \sum_{i=1}^K \frac{m_i}{2} v_i^{\mathrm{T}} v_i +  \frac{1}{2} \omega_i^{\mathrm{T}} R_i  \mathbb{I}_i R_i^{\mathrm{T}} \omega_i  - U(x_i, R_i) \text{.}
\]
Note that $l(x_i, v_i, R_i, \omega_i) $ is shorthand notation for $l(x_1, v_1, R_1, \omega_1, \cdots, x_K, v_K, R_K, \omega_K)$.  We will use this shorthand notation elsewhere to simplify the expressions. The path that the stochastic mechanical system takes on the time-interval $[a,b]$ is one that extremizes the HP action:
\[
s = \int_a^b \left[ l(x_i, v_i, R_i, \omega_i) dt + \sum_{q=1}^m \gamma_q(x_i, R_i) \circ dW_q + \left\langle p_i, \frac{d x_i}{dt} - v_i \right\rangle dt + \left\langle \widehat{\pi}_{i},  \frac{d R_i}{dt} R_i^{-1} - \widehat{\omega}_i  \right\rangle dt \right] \, 
\]
for arbitrary variations with fixed endpoints: $(x_i(a), R_i(a))$ and $(x_i(b), R_i(b))$.  The corresponding SDEs of motion are given by:
\begin{align*}
d x_i  = v_i dt~~~ & \text{(reconstruction equation),} \\
d p_i = - U_{x_i} dt + \sum_{q=1}^m (\gamma_q)_{x_i} \circ dW_q(t, \omega)~~~ & \text{(stochastic EL equations),} \\
p_i = m_i v_i ~~~ & \text{(Legendre transform), } \\
d R_i = \widehat{\omega_i} R_i dt~~~ & \text{(reconstruction equation), } \\
d\pi_i = - U_{R_i} dt + \sum_{q=1}^m (\gamma_q)_{R_i} \circ dW_q~~~ & \text{(stochastic LP equations),} \\
\pi_i  = R_i  \mathbb{I}_i R_i^{\mathrm{T}} \omega_i~~~ & \text{(reduced Legendre transform).}
\end{align*}
for $i=1,...,K$.  The terms $U_{x_i}$ and $U_{R_i}$  are defined in terms of the inner product on $\mathbb{R}^3$ as,
 \begin{align*}
 U_{x_i}^{\mathrm{T}} y &= \left\langle \frac{\partial U}{\partial x_i}, y \right\rangle  = \partial_{x_i} U(x_i, R_i) \cdot y \text{,} \\
U_{R_i}^{\mathrm{T}} y &= \left\langle \frac{\partial U}{\partial R_i} R_i^{\mathrm{T}}, \widehat{y} \right\rangle =   \partial_{R_i} U(x_i, R_i) \cdot \widehat{y} R_{i} \text{,}
 \end{align*}
where $\partial_{R_i} U:  \operatorname{SO}(3) \to T_{R_i}^*\operatorname{SO}(3)$, and  $\partial_{x_i} U: \mathbb{R}^3 \to T_{x_i}^* \mathbb{R}^3$ as defined below.
Adding dissipation so that the Gibbs distribution is invariant under the stochastic process defined by the above SDE with dissipative drift, yields Langevin-type equations for rigid-body systems (see, e.g., \citep{BoOw2007c}).

\paragraph{Stochastic Variational Euler integrator}

For the discrete description, the variational Euler integrator provided earlier is implemented.   Let $B_q^k \sim \mathcal{N}(0,h)$ be normally distributed random variables for $q=1,...,m$ and $k=0,...,N-1$.   The action sum is given by 
\begin{align*}
s_d = \sum_{k=0}^{N-1}  & \left[ \ell \left(x_i^k, v_i^k,  R_i^k, \omega_i^k \right)  h
+ \left\langle p_i^{k+1}, (x_i^{k+1}-x_i^k)/h - v_i^{k+1} \right\rangle  h \right] \\
&+  \left[ \left\langle  \widehat{\pi_i^{k+1}}, \tau^{-1}(R_i^{k+1} (R_i^k)^\mathrm{T} )/h -\widehat{\omega_i^{k+1}} \right\rangle h \right] + \sum_{q=1}^m \gamma_q(x_i^k, R_i^k) B_q^k \text{.}
\end{align*}
Stationarity of this action sum implies the following discrete scheme:
\begin{align*}
x_i^{k+1} &= x_i^k + h v_i^{k+1} ~~~ \\
p_i^{k+1} &= p_i^k - h U_{x_i}(x_i^k, R_i^k) +  \sum_{q=1}^m (\gamma_q)_{x_i}(x_i^k, R_i^k) B_q^k~~~ \\
p_i^k &= m v_i^k ~~~ \\
R_i^{k+1} &= \tau \left(\widehat{\omega_i^{k+1}} h \right) R_i^{k} ~~~  \\
\left(d \tau^{-1}_{h  \omega_i^{k+1}} \right)^* \pi_i^{k+1} &=  \left(d \tau^{-1}_{h \omega_i^k} \right)^* \pi_i^k -  h U_{R_i}(x_i^k, R_i^k) +  \sum_{q=1}^m (\gamma_q)_{R_i}(x_i^k, R_i^k) B_q^k \\
\pi_i^k &= R_i^k \mathbb{I}_i  (R_i^k)^{\mathrm{T}} \omega_i^k  ~~~ 
\end{align*}
for $i=1,...,K$.  Assuming the Legendre transforms are invertible, this integrator has the attractive property that the translational and rotational configuration updates, and the translational momentum update, are explicit.  One only has to perform an implicit solve for the discrete Lie-Poisson part.  Even that computation is straightforward since the torque due to the potential is only a function of the orientation and position at the previous time-step.

\section*{Acknowledgements}We wish to thank Andreu Lazaro, Jerry Marsden, Juan-Pablo Ortega and Katie Whitehead for stimulating discussions.

\bibliographystyle{plain}
\bibliography{nawaf}

 \end{document}